\documentclass[12pt]{amsart}

\usepackage{amssymb}

\newtheorem{thm}{Theorem}%[section]

\newtheorem{lem}[thm]{Lemma}

\newtheorem{prop}[thm]{Proposition}

\newtheorem{conj}[thm]{Conjecture}
\theoremstyle{definition}
\newtheorem{defn}[thm]{Definition}

\newtheorem{say}[thm]{}
\newtheorem{exmp}[thm]{Example}

\newtheorem{rem}[thm]{Remark}          
\newtheorem{ack}{Acknowledgments}

\newtheorem{defn-thm}[thm]{Definition--Theorem}  %!!!!!!!!!!!!!!!!!!!!!!!!
\theoremstyle{remark}
\newtheorem{claim}[thm]{Claim}

\renewcommand{\o}[0]{{\mathcal O}} 
\newcommand{\z}[0]{{\mathbb Z}}

\renewcommand{\a}[0]{{\mathbb A}}

\newcommand{\p}[0]{{\mathbb P}}
\newcommand{\f}[0]{{\mathbb F}}
\newcommand{\q}[0]{{\mathbb Q}}
\newcommand{\map}[0]{\dasharrow}
\newcommand{\qtq}[1]{\quad\mbox{#1}\quad}
\newcommand{\spec}[0]{\operatorname{Spec}}

\newcommand{\im}[0]{\operatorname{im}}

\newcommand{\Hom}[0]{\operatorname{Hom}}

\newcommand{\chr}[0]{\operatorname{char}}

\newcommand{\hilb}[0]{\operatorname{Hilb}}

\newcommand{\onto}[0]{\twoheadrightarrow}

\newcommand{\ch}[0]{\operatorname{CH}}

\def\into{\DOTSB\lhook\joinrel\rightarrow}

\begin{document}
\bibliographystyle{amsplain}

\title{Specialization of zero cycles}
\author{J\'anos Koll\'ar}

\maketitle

Let $X$ be a proper scheme over a field $K$.
There are two ways of organizing the points of $X$
into equivalence classes using rational curves.
  One is the notion of {\it R-equivalence}
introduced by \cite{manin}. Two points
$x_1,x_2\in X(K)$ are called  directly R-equivalent
if there is a morphism $p:\p^1\to X$ with $p(0{:}1)=x_1$ and
$p(1{:}0)=x_2$.  This generates an equivalence relation called 
 R-equivalence.
The set of R-equivalence classes $X(K)/R$
forms a set.

Closely related to it is {\it rational equivalence}.
Here we allow pairs of morphisms $h:C\to \p^1$ and $p:C\to X$
and declare $p_*(h^{-1}(0{:}1))$ and $p_*(h^{-1}(1{:}0))$
to be rationally equivalent.
Rational equivalence is sometimes  hard to
see geometrically.
The set of rational equivalence classes forms a group $\ch_0(X)$.

 Let $S$ be the spectrum of a local Dedekind ring with
residue field $k$ and quotient field $K$.
Let $X_S\to \spec S$ be a proper morphism,
  $x_K$  a closed point of $X_K$ and $x_S$ its closure in $X_S$.
Then $x_S\cap X_k$ is a zero cycle on $X_k$. This defines a
specialization map on zero cycles.
It is easy to prove (see, for instance, \cite[2.3]{fulton}),
that this descends to   {\it specialization maps}
$$
\ch_0(X_K)\to \ch_0(X_k)\qtq{and} X_K(K)/R\to X_k(k)/R.
$$
In general there is very little that one can say about these maps.

If $X_S\to \spec S$ is a family of curves,
then the specialization maps are neither surjective
nor injective.

The only reasonable case when surjectivity holds is
when $X_S$ is smooth over $S$ and $S$ is Henselian.
If this is assumed then injectivity holds only if 
  $X_S \cong \p^1_S$.

The aim of this paper is to study 
higher dimensional cases where the specialization map
is an isomorphism. 
The correct higher dimensional analogs of rational
curves are the 
{\it separably rationally connected} or SRC varieties.
See \cite{Ko01} for an introduction to their theory
and \cite[IV.3]{rcbook} for a more detailed treatment.
There are many equivalent conditions defining this notion.
The definition given  below essentially  says that
two general points can be connected by a rational curve.
In positive characteristic we also have to be mindful
of some inseparability problems.

\begin{defn}\label{SRC-defn}  A smooth, 
proper variety $X$ is called {\it separably rationally connected}
  or {\it SRC}, if there is a variety $U$ and a
morphism $F:U\times\p^1\to X$ such that the induced map
$$
F(\_\ ,(0{:}1))\times F(\_\ ,(1{:}0)) : U\to X\times X
$$
is dominant and separable.
\end{defn}

The main result is the following:

\begin{thm} \label{main.thm}
Let $S$ be the spectrum of a local, Henselian, Dedekind ring with
residue field $k$ and quotient field $K$.
Let $X_S\to \spec S$ be a smooth proper morphism.
Assume that $X_k$ is SRC. Then
\begin{enumerate}
\item The specialization map on R-equivalence
$$
X_K(K)/R\to X_k(k)/R\qtq{is an isomorphism of sets.}
$$
\item If $k$ is perfect, the specialization map on the 
Chow group of zero cycles
$$
\ch_0(X_K)\to \ch_0(X_k)\qtq{is an isomorphism of groups.}
$$

\end{enumerate}
\end{thm}

If $k$ is not perfect, we sill get that
$\ch_0(X_K)\to \ch_0(X_k)$ is an isomorphism modulo $p$-torsion
where $p=\chr k$.
I do not know any example where we do not have an isomorphism.

The second part of (\ref{main.thm}) has been known in case $k$ is finite.
For cubic hypersurfaces this is done in \cite{madore},
the general case is treated in \cite{ko-sz}. In this case
both sides are trivial.
The first part has been known in case $k$ is finite
and sufficiently large (depending on the dimension and degree of $X$
under some projective embedding),
 see \cite{ko-sz}.
In all of these cases the idea is to prove that
the triviality of the Chow group or of  R-equivalence
over a finite field is shown by some ``very nice'' maps
$p:\p_k^1\to X_k$, and these can be lifted to $X_K$.

In general, we have to deal with the lifting of
arbitrary maps $\p_k^1\to X_k$. This is accomplished
using an observation of \cite[2.4]{ghs}.

A theory for lifting of families of torsors over finite group schemes
was developed in the papers \cite{mb01, mb02}.  
Roughly speaking, this method implies   (\ref{main.thm}.1)
for varieties birational to a quotient  $\a^n/G$.

\medskip

As an application of the method of the proof of
(\ref{main.thm}) we obtain 
two further results. The first is 
an upper semi continuity
statement for the number of R-equivalence classes
in families of varieties over local fields:

\begin{thm}\label{semicont.thm} Let $K$ be a local field and
$f:X\to Y$ a smooth, projective $K$-morphism
whose fibers are SRC. Then
$$
Y(K)\ni y\mapsto |X_y(K)/R| \qtq{is upper semi continuous}
$$
in the $p$-adic topology (\ref{loc.field.defn}).
\end{thm}

\begin{rem} It is quite likely that $y\mapsto |X_y(K)/R|$
is actually a continuous function
(that is, locally constant). All varieties with a
given reduction form an open set, thus
the continuity follows from (\ref{main.thm}) if every fiber $X_y$
has a smooth SRC reduction.

It is also likely that the analogously defined function
$y\mapsto |\ch^0_0(X_y)|$ is also continuous.
\end{rem}

The second applications asserts that R-equivalence and
direct R-equivalence coincide over large fields (\ref{large.def}).
The precise assertion  (\ref{large.fld.thm}) is actually much stronger.

\medskip

The third part of the paper    uses specialization to  singular
varieties to get two types of examples of quartic 
hypersurfaces:

\begin{exmp}\label{nonfg.exmp} 
For every $n\geq 5$
 there is a smooth quartic hypersurface
$H^n\subset \p^{n+1}$ over $\q(t)$ 
such that $H^n$ is  unirational (over $\q(t)$) and it has
 infinitely many
R-equivalence classes. 
\end{exmp}

By a theorem of Springer, if a smooth quadric hypersurface
defined over $K$ has a point
in an odd degree extension of $K$, then it has a point in $K$
itself.  Similarly, it was conjectured by 
Cassels and Swinnerton-Dyer
that if a smooth cubic
hypersurface of dimension $\geq 2$  over $K$ has a point
in an extension of $K$ whose degree is not divisible by 3,
 then it has a point in $K$
itself.  Many cases of  this have been proved in \cite{coray}.

The next example shows that no similar result holds for quartics.

\begin{exmp}\label{deg1.zeroc.exmp} 
 For every odd $d$ and $n\gg 1$ there is a smooth quartic
$H_d^n\subset \p^{n+1}$ over $K=\q(t)$ such that 
$H_d^n$ has a point
in a  field extension  of degree $d$ but
it does not have a point in any
 field extension  of smaller odd degree.
\end{exmp}

\section{Deformations of combs}

Let $X_S\to \spec S$ be a smooth morphism as in
(\ref{main.thm}) and 
$x_K,y_K\in X_K(K)$ points with specializations $x_k,y_k\in X_k(k)$.
Let $C_k\subset X_k$ be a
rational curve showing that $x_k,y_k$  are R-equivalent.
(\ref{main.thm}.1) follows if every such $C_k$
can be lifted to a rational curve $C_S\subset X_S$
passing through $x_K,y_K$.
This is a deformation theory problem
with known obstructions in the first cohomology
of the normal bundle of $C_k$ 
twisted by $\o_{C_k}(-x_k-y_k)$ (if $C_k$ is smooth).  
We run into difficulties if this first cohomology group
is not zero.

Following the method of \cite{kmm2},
we try  to deal with this problem by attaching
auxiliary curves $A_i$ to $C_k$ and deform the
resulting reducible curve $C^*_k:=C_k\cup A_1\cup\cdots\cup A_m$.

It was observed in \cite{kmm2} that for
suitable $A_i$ and  large $m$ the
deformation theory of $C^*_k\subset X$ is  better 
behaved than the deformation theory of $C_k$ itself.
\cite{kmm2} concentrated on the deformations
of the morphism $C^*_k\to X$, in which case the obstructions
never vanish.

Recently,  \cite{ghs} proved that if we look at the Hilbert scheme
instead, then obstructions vanish for suitable choices 
of  the $A_i$.
For arithmetic applications this is a crucial
improvement.

\begin{defn} \label{defcomb} 
Let $C$ be a geometrically reduced projective curve
over a field $k$. 
A {\it comb} over $C$ with {\it $n$-teeth} is a 
reduced projective curve $C\cup A_1\cup\cdots \cup A_n$
having $n$ more  irreducible components over $\bar k$.
$C$ is called the {\it handle}.
The other $n$ components, $A_1,\dots, A_n$ are smooth rational curves, 
 disjoint from each other and 
intersect  $C$ transversally 
in $n$ distinct smooth points. The curves $A_1,\dots, A_n$ may not be 
individually defined
over $k$. 
A comb can be pictured as below:
$$
\begin{array}{c}
\begin{picture}(100,100)(40,-70)
\put(0,0){\line(1,0){180}}
\put(20,10){\line(0,-1){60}}
\put(40,10){\line(0,-1){60}}
\put(130,10){\line(0,-1){60}}
\put(160,10){\line(0,-1){60}}

\put(15,-60){$A_1$}
\put(35,-60){$A_2$}
\put(125,-60){$A_{n-1}$}
\put(155,-60){$A_n$}

\put(-25,-5){$C$}

\put(70,-30){$\cdots\cdots$}

\end{picture}\\
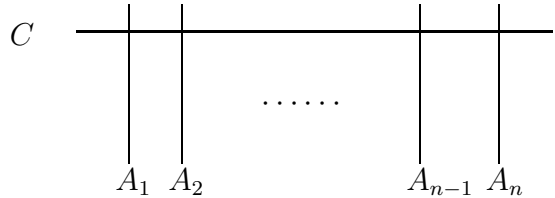
\mbox{Comb with $n$-teeth}
\end{array}
$$
\end{defn}

\begin{say}[Construction of combs]\label{const.comb}

Let $X$ be a smooth projective variety over
an algebraically closed field $\bar k$.
Let $C\subset X$ be a reduced local complete intersection curve
with ideal sheaf $I_C$.
(In our applications, $C$ will have only  nodes.)
Let $N_C:=Hom(I_C/I_C^2,\o_C)$ denote the normal bundle.

The tangent space of the Hilbert scheme of $X$ at $[C]$
 is $H^0(C,N_C)$ and the obstructions
lie in $H^1(C,N_C)$. Our aim is to create a comb $C^*$ with handle $C$
such that $H^1(C^*,N_{C^*})$.

To do this, assume that we have a collection of smooth rational curves
$\{A_w:w\in W\}$
such  that the following 3 conditions are satisfied:
\begin{enumerate}
\item The normal bundle $N_{A_w}$ is  
semi positive for every $w$. (That is, it is
a direct sum of  
line bundles of nonnegative degrees.)
\item Every $A_w$ intersects $C$ in a single point,
and the intersection is transverse. Let  this point  be $p=p_w$,
 $L_{p,w}\subset T_{p,X}$  the
 tangent line of $A_w$ at $p$ and  $\bar L_{p,w}$
 its image in the fiber $N_{p,C}$ of $N_C$ over $p$.
\item There is a dense set of smooth points $p\in C$
such that the lines 
$$
\bar L_{p,w}: p\in A_w\ \mbox{span}\ N_{p,C}\ \mbox{(as a vector space).}
$$
\end{enumerate}
\end{say}

\begin{rem}  For an algebraically closed field, the
collection of smooth rational curves $\{A_w:w\in W\}$
can be chosen to be an algebraic family. In the applications
to nonclosed fields, however,  the collection $\{A_w:w\in W\}$
will not be an algebraic family, rather a  subset of an algebraic
family determined by certain arithmetic conditions, see 
(\ref{nonclosed.say}).

In practice, the first condition is usually
easy to satisfy. If $F:W\times \p^1\to X$
is a dominant and separable morphism then 
$F^*T_X$ is semi positive on the general $\p^1_w$
and so is the relative normal sheaf $F^*T_X/T_{\p^1}$.
If $F$ is an embedding on the general $\p^1_w$,
 we can pass to an open set $W^0\subset W$
to ensure that (\ref{const.comb}.1) holds. 

(\ref{const.comb}.2) is also easy to achieve 
 if $\dim X\geq 3$, see \cite[II.3.]{rcbook}.
We need rational curves  through a point with ample normal bundle.
These exist an any SRC variety.

The key condition is (\ref{const.comb}.3). This can be satisfied in
two important cases.

First, if $X$ is SRC. Indeed, in this case for every
$p\in X$ and every tangent direction
$v\in T_pX$ there is a rational curve
$f:\p^1\to X$ such that $f(0{:}1)=p$ and it has tangent direction $v$ there.

Second, assume that there is a morphism $h:X\to B$ to a smooth curve $B$
whose general fibers are SRC. Let $C\subset X$ be a curve such that
$h:C\to B$ is separable.  At a general point $p\in C$
the fiber of the normal bundle $N_{p,C}$ can be identified with
$T_pX_p$ where $X_p$ denotes the fiber through $p$.
The SRC case  now gives enough rational curves
in $X_p$ to satisfy (3).
\end{rem}

Assume now that we have $X,C$ and the curves 
$\{A_w:w\in W\}$
as in (\ref{const.comb}.1--3). 
Pick $w_i\in W, i=1,\dots,m$ such that the  $A_{w_i}$
intersect $C$ in {\em distinct}  smooth points.
Then the union 
$$
C^*=C(w_1,\dots,w_m):=C\cup A_{w_1}\cup \cdots \cup A_{w_m}
$$
is a reduced curve which is a local complete intersection along $C$.

\begin{lem}\cite{ghs}\label{hgs-lem}
Let $X$ be a smooth projective variety,
$C\subset X$ a geometrically reduced local complete intersection curve
and $\{A_w:w\in W\}$ a collection of curves
satisfying the conditions (\ref{const.comb}.1--3).
Let $M$ be any line bundle on $C$.

Then there are $w_1,\dots,w_m\in W$ such that
 the resulting curve $C^*=C(w_1,\dots,w_m)$
satisfies the following three  conditions:
\begin{enumerate}
\item $N_{C^*}|_C$ is generated by global sections.
\item $H^1(C,M\otimes N_{C^*}|_C)=0$.
\item The above two conditions  hold for any
other $C(w_1,\dots,w_m,\dots,w_{m+n})$
obtained by adding more curves to $C^*$.
\end{enumerate}
\end{lem}

Proof. By Serre duality, 
$$
H^1(C,M\otimes N_{C^*}|_C)\qtq{is dual to}
\Hom(M \otimes \omega_C^{-1}, I_{C^*}/I^2_{C^*}|_C).
$$
The key point is to understand how $I_{C^*}/I^2_{C^*}|_C$
changes if we add a new curve $A_{w_{m+1}}$
with intersection point $p=p_{m+1}$.
Set $C^*(w_{m+1})=C(w_1,\dots,w_{m+1})$.
Let  $L(p)\subset N_{p,C}$ be the
tangent line of $A_{w_{m+1}}$ with dual map
$q_p:N_{p,C}^*\to L(p)^*$. Let further
$r_p:I_{C^*}/I^2_{C^*}|_C\to N_{p,C}^*$ denote the restriction.
We have an exact sequence
$$
0\to I_{C^*(w_{m+1})}/I^2_{C^*(w_{m+1})}|_C\stackrel{r_p}{\to}  
I_{C^*}/I^2_{C^*}|_C
\stackrel{q_p\circ r_p}{\longrightarrow} L(p)^*\to 0.
$$
Pick any  
$$
\phi\in  \Hom(M \otimes \omega_C^{-1}, I_{C^*}/I^2_{C^*}|_C).
$$
There is an open  set of points 
$p\in C$ such that $\phi$ has rank one at $p$.
Thus the composition 
$$
r_p\circ \phi: 
M \otimes \omega_C^{-1}\otimes k(p)\to N_{p,C}^*
$$
is an injection. By assumption (\ref{const.comb}.3), there is a curve
$A_{w_{m+1}}$ such that the induced
map
$$
q_p\circ r_p\circ \phi: 
M \otimes \omega_C^{-1}\otimes k(p)\to 
N_{p,C}^*\to L(p)^*
$$
is nonzero. Thus 
$$
\phi\notin \Hom(M \otimes \omega_C^{-1},
I_{C^*(w_{m+1})}/I^2_{C^*(w_{m+1})}|_C).
$$
Hence by adding suitable curves  $A_{w_s}$
we eventually achieve (\ref{hgs-lem}.2)
and the vanishing still holds by adding more points.

Global generation of $N_{C^*}|_C$ can be guaranteed
by the vanishing of $H^1(C,\o_C(-P)\otimes N_{C^*}|_C)$
for all $P$, which in turn follows from the
vanishing of the single cohomology group
$H^1(C,L^{-1}\otimes N_{C^*}|_C)$
where $L$ is any globally generated line bundle.\qed

\begin{say}[Dealing with extra intersections]\label{get.rid.int}

In all of our  cases one can arrange that the $A_{w_i}$ are disjoint
from each other. Then $C^*$ is a comb in $X$ and this is what we want.

In general we may have 
 two curves, say $A_i$ and $A_j$, intersecting  at  a point
$x\in X\setminus C$. 
One can get around such extra intersections  with 
one of the following tricks.

First, we can
replace $X$ with a smooth but nonseparated scheme $X'$
where the point $x$ is replaced by 2 points $x_i$ and $x_j$.
We can lift $A_i$ to $A'_i$ going through $x_i$ and $A_j$ to
$A'_j$ going through $x_j$. This way we remove an intersection point.
The local deformation theory of the Hilbert scheme
is the same on nonseparated schemes, and this is all we need.
The global theory is more problematic but we do not use it.

Second, we can
replace $X$ with $X'=X\times \p^1$. We  then replace
the curves $A_w$ by graphs of suitable
isomorphisms in $A_w\times \p^1$. We also have to add the
trivial rational curve $\{p\}\times \p^1$ to our collection $\{A_w\}$.

If $X$ itself is of the form $Y\times \p^n$ for some $n\geq 1$
then we can move curves in $X$ itself, so  a further
product with $\p^1$ is not needed.

Both of these methods are suitable for our current purposes.
\end{say}

\begin{prop}\cite{ghs}\label{hgs-prop}
Let $X$ be a smooth projective variety,
$C\subset X$ a geometricaly reduced local complete intersection curve
and $\{A_w:w\in W\}$ a collection of curves
satisfying the conditions (\ref{const.comb}.1--3).
Let $M$ be any line bundle on $C$.

Then there are $w_1,\dots,w_m\in W$ such that
(after possibly passing to a scheme $X'$ as in (\ref{get.rid.int})) 
the resulting comb $C^*=C(w_1,\dots,w_m)$
satisfies the following three  conditions:
\begin{enumerate}
\item $N_{C^*}$ is generated by global sections.
\item $H^1(C^*,M^*\otimes N_{C^*})=0$, where $M^*$ is the unique extension
of $M$ to $C^*$ which has
degree 0 on the other components.
\item The above two conditions  hold for any
other $C(w_1,\dots,w_m,\dots,w_{m+n})$
obtained by adding more curves to $C^*$.
\end{enumerate}
\end{prop}

Proof. Everything follows from (\ref{hgs-lem}) and
the exact sequences
$$
0\to \sum_{i} N_{C^*}|_{A_{w_i}}\otimes \o(-p_i)
\to N_{C^*}
\to N_{C^*}|_C\to 0
$$
and
$$
0\to N_{A_{w_i}}\to N_{C^*}|_{A_{w_i}} \to T_{p_i,C}\to 0.\qed
$$

\begin{say}[Nonclosed fields]\label{nonclosed.say}

In our applications we work over a field $k$
which is not algebraically closed
and we want to obtain a  curve $C^*$
which is defined over $k$.
Assume that we have a solution 
$C^*=C(w_1,\dots,w_m)$ of (\ref{hgs-lem})
over $\bar k$. Let
$A_{w_1},\dots,A_{w_s}$ denote all the conjugates 
of the curves $A_{w_1},\dots,A_{w_m}$. Then
$C^{**}:=C(w_1,\dots,w_s)$
is defined over $k$. Moreover,
(\ref{hgs-lem}) applies if the curves
$A_{w_1},\dots,A_{w_s}$ intersect $C$ in
{\em distinct} points.

This is a somewhat  troublesome
condition which is arithmetic in nature. The
problem is essentially the following:
\begin{enumerate}
\item[{}]
Given $C\subset X$, find points $p\in C$ and  rational curves
$p\in A\subset X$ defined over $k(p)$. We also want $k(p)/k$
to be separable.
\end{enumerate}
It is clearly enough to consider this problem for irreducible curves.

\end{say}

Let $C$ be an irreducible curve, $C_{gen}\in C$  the generic point and 
$v_1,\dots,v_n$  a basis of the fiber of  $N_C$ over $C_{gen}$.
Let $H(C_{gen},v_i)$ be an irreducible component of
the Hilbert scheme of smooth rational curves in $X$ with ample normal 
bundle intersecting $C$ only at  $C_{gen}$ with tangent direction $v_i$.
Each $H(C_{gen},v_i)$ is a smooth variety over $k(C_{gen})=k(C)$
which can be extended to a $k$-variety
$$
\tau_i: H(C,v_i)\to C
$$
which parametrizes smooth rational curves in $X$ with ample normal 
bundle intersecting $C$ in  a single point and whose tangent vector
at that point is the corresponding specialization of $v_i$.

We need to find points $p_j\in C$ such that
each $\tau_i^{-1}(p_j)$ has a $k(p_j)$-point.

Since $\tau_i$ is generically smooth,
 we can take curve sections $B_i\subset H(C,v_i)$
such that $\tau_i:B_i\to C$ is dominant and separable.
Let 
$$
B\subset B_1\times_C \cdots \times_C B_n
$$
be an irreducible component with projection $\tau:B\to C$.

\begin{lem} Let $k$ be a field and $f:B\to C$  a dominant morphism of 
smooth 
curves over $k$. Then there is a dense set of points $q_i\in B$
such that $k(q_i)=k(f(q_i))$ and $k(q_i)/k$
is separable.
\end{lem}

Proof. Pick a point $c\in C$ such that $f$ is finite over $c$.
Let $h:f^{-1}(c)\to \a^1$ be an embedding and extend it to a separable
rational  function $h$ which is not constant on any
geometric irreducible component of $B$. $(f,h)$
gives a morphism $C\to C'\subset B\times \p^1$
which is birational onto its image. 

Let $u\in \p^1$ be any point such that the second projection 
$\pi_2:C'\to \p^1$
is smooth over $u$. Then any $q\in \pi_2^{-1}(u)$ has the required property.
\qed 
\medskip

Putting everything together, we obtain 
the following results about the existence of ``good'' combs
over arbitrary fields.

\begin{thm} \label{SRC.noobs.thm}
Let $X$ be a smooth, proper, SRC variety over a field $k$.
Let $C$ be a 
geometrically reduced local complete intersection curve,
$S\subset C$ a finite set of smooth points.
Let $g:C\to X$ be a morphism and $G:C\into X\times \p^n$
an embedding lifting $g$.

Then there is a comb $C^*\subset X\times \p^n$
 defined over $k$ with handle $G(C)$
such that
\begin{enumerate}
\item $N_{C^*}$ is generated by global sections.
\item $H^1(C^*, \o_{C^*}(-G(S))\otimes N_{C^*})=0$.
\item $C^*$ is smooth at $G(S)$.\qed
\end{enumerate}
\end{thm}

\begin{thm}\label{SRCfiber.noobs.thm}
 Let $X$ be a smooth, proper,  variety over a field $k$
and $f:X\to B$ a morphism to a smooth curve whose general fibers are SRC.
Let $C$ be a 
geometrically reduced local complete intersection curve,
$S\subset C$ a finite set of smooth points.
Let $g:C\to X$ be a morphism 
such that $f\circ g:C\to B$ is separable
and $G:C\into X\times \p^n$
an embedding lifting $g$.

Then  there is a comb $C^*\subset X\times \p^n$ 
defined over $k$ with handle $G(C)$
whose teeth are contained in fibers of $f$
such that
\begin{enumerate}
\item $N_{C^*}$ is generated by global sections.
\item $H^1(C^*, \o_{C^*}(-G(S))\otimes N_{C^*})=0$.
\item $C^*$ is smooth at $G(S)$.\qed
\end{enumerate}
\end{thm}

\section{Proof of the main theorem}

The specialization map on zero cycles
is surjective since every point in $X_k$ has a lifting
to $X_K$ by the Hensel property.
Thus it is enough to prove the following 2 assertions: 
\begin{enumerate}
\item If the specializations of 2 points are directly R-equivalent,
then the 2 points are also R-equivalent, and
\item If the specializations of 2 zero cycles  are rationally  equivalent,
then the 2 zero cycles  are also  rationally equivalent.
\end{enumerate}

\begin{say}[Proof of (\ref{main.thm}.1)]\label{main.2.pf}

Let $u_K,v_K\in X_K(K)$ be 2 points with
specializations $u_k,v_k$. If $u_k,v_k$ are directly R-equivalent, then
there is a morphism
 $f:\p_k^1\to X_k$ 
such that $f(0{:}1)=u_k$ and $f(1{:}0)=v_k$.
(If $u_k=v_k$ then we choose the constant map.)

Consider the diagonal map
$$
F:\p^1\to X\times \p^1, \quad  p\mapsto (f(p),p).
$$
Its image is a smooth rational curve $C\subset X\times \p^1$
connecting $u'_k:=(u_k,(0{:}1))$ and $v'_k:=(v_k,(1{:}0))$.

By (\ref{SRC.noobs.thm}) there is a comb $C^*$ defined over $k$
such that $C^*$ is smooth at $u'_k,v'_k$ and 
$H^1(C^*,\o_{C^*}(-u'_k-v'_k)\otimes N_{C^*})=0$.

Mostly for ease of reference, we replace
$X\times \p^1$ by its blow up $Z_S:=B(X_S\times \p^1)$ along
the 2 sections $(u_S,(0{:}1))$ and $(v_S,(1{:}0))$.
Let $E_u,E_v\subset Z$ denote the exceptional divisors.
Let $H_S\to \spec S$ denote the relative Hilbert scheme of 
one dimensional subschemes of $Z_S$.

Let $[C^*]\in H_S$ be the point corresponding to 
(the birational transform of) $C^*$.
By the theory of Hilbert schemes (cf.\  \cite[I.2.15]{rcbook}
$H^1(C^*,\o_{C^*}(-u_k-v_k)\otimes N_{C^*})=0$
implies that 
$H_S\to \spec S$ is smooth at $[C^*]$.
By the Hensel property, there is  a morphism
$\sigma:\spec S\to H_S$ such that $\sigma(\spec k)=[C^*]$.
Thus $\sigma(\spec K)$ corresponds to a reduced genus zero curve $C_K$,
defined over $K$
which has intersection number 1 with $E_u$ and $E_v$.
This implies that  $C_K\cap E_u$ and $C_K\cap E_v$ are $K$-points of $C_K$.
The projection of $C_K$ to $X_K$ therefore shows that
$u_K$ and $v_K$ are directly R-equivalent.\qed
\end{say}

\begin{say}[Proof of (\ref{main.thm}.2)]\label{pf.of.main.1}

This is done in 2 steps.
\begin{enumerate}
\item 
First we prove that if 2 zero cycles $Z^1_K,Z^2_K\in X_K$
have the same specializations, then they are
rationally equivalent.
\item
Then we prove that if a zero cycle $Z_k\in X_k$
is rationally equivalent to 0, then it has a lifting
$Z_K$ which is also  rationally equivalent to 0.
\end{enumerate}

In order to see (\ref{pf.of.main.1}.1),
write 
$$
Z_K:=Z^1_K-Z^2_K=\sum_j m_j[P_j]
$$
 as a sum of irreducible zero cycles.
For a point $p\in X_k$ let 
$$
Z_K(p)=\sum^{(p)}_jm_j[P_j]
$$
 be the sum of those points
whose specialization is $p$. 
It is  enough to prove that
$Z_K(p)$ is rationally equivalent to 0 for every $p$.

The next  argument follows a  suggestion of  
Colliot-Th\'el\`ene.

The specialization of $Z_K(p)$ is 0, thus
$\sum^{(p)}_jm_j\deg [P_j]=0$.
Let $P\in X_K$ be a lifting of $p$. It is then sufficient 
to  prove  that
$$
[P_j]-\frac{\deg [P_j]}{\deg [P]}[P]
$$
is rationally equivalent to 0 for every $P_j$ specializing to $p$.

Take a base field extension from $K$ to $K(P_j)$.
From (\ref{hensel.ext.lem}) we see that $K(P)\subset K(P_j)$
thus both $P_j$ and $P$ become $K(P_j)$-points
which specialize to the same point of $X_{k(p)}$.
By (\ref{main.thm}.1) these two points are
R-equivalent, thus $[P_j]-[P]\in \ch_0(X_{K(P_j)})$ is zero.
Its push forward to $X_K$ is 
$[P_j]-(\deg [P_j]/\deg [P])[P]$, which is thus also zero in $\ch_0(X_K)$.
\medskip

The proof of (\ref{pf.of.main.1}.2) 
uses (\ref{SRCfiber.noobs.thm}). 

By assumption, there is a smooth curve and  morphisms
$f:C\to \p^1$ and $g:C\to X_k$ such that
$Z_k=g_*(f^{-1}(0{:}1))-g_*(f^{-1}(1{:}0))$.

If $f$ is not separable, then we can factor it as
a separable map $f_1:C\to \p^1$ composed with a purely inseparable map 
$\phi:\p^1\to\p^1$
(this is always possible for curves over a perfect field). This shows that 
$$
Z_k = \deg \phi \cdot [g_*(f_1^{-1}(0{:}1))-g_*(f_1^{-1}(1{:}0))],
$$
thus it is sufficient to consider the case when $f$ is separable.

Choose an embedding $j:C\into \p^m$ and 
apply (\ref{SRCfiber.noobs.thm}) to
the morphism 
$$
X\times \p^m\times \p^1\to \p^1\qtq{with} C':= \im (g\times j\times f).
$$
There is a comb $C^*$ with handle $C'$
which is unobstructed in  the Hilbert scheme
as in (\ref{main.2.pf}). 
This gives a curve $C^*_K\subset X_K\times \p^m\times \p^1$
with projections $F:C^*_K\to \p^1_K$ and $G:C^*_K\to X_K$.
By our construction
the degree of $F:C^*\to \p^1$ is the same as the degree of $f:C\to \p^1$,
thus
$$
Z_K:=G_*(F^{-1}(0{:}1))-G_*(F^{-1}(1{:}0))
$$
is a zero cycle on $X_K$ which is
rationally equivalent to zero and whose specialization is $Z_k$.\qed

\begin{lem}\label{hensel.ext.lem}
 Let $(R,m)$ be a Henselian local ring.
\begin{enumerate}
\item There is a one--to--one correspondence
between finite \'etale $R$-algebras and
finite separable $R/m$-algebras.
\item Let $S$ be a finite, flat  $R$-algebra such that $S/\sqrt{mS}$
is separable over $R/m$.
Then there is a finite \'etale $R$-algebra $S'\subset S$
such that $S/\sqrt{mS}\cong S'/\sqrt{mS'}$.
\end{enumerate}
\end{lem}

Proof. If $(R,m)$ is complete, this is essentially
the theory of Witt rings as explained in
\cite[Chap.II]{serre-loc}.

In the general case,
for the first part see \cite[I.4.4]{milne}.
To see the second part, let $S'$ be the finite \'etale $R$-algebra
corresponding to $S/\sqrt{mS}$. Consider  $S'\otimes_R S$.
This is a finite, \'etale $S$-algebra and $S$ is also Henselian
by \cite[I.4.3]{milne}. $S/\sqrt{mS}$ is a direct summand of
$S/\sqrt{mS}\otimes_{R/m}S/\sqrt{mS}$,
thus by applying the first part to $S'\otimes_R S$ over $S$
we obtain that $S$ is a direct summand of $S'\otimes_R S$.
This gives the required embedding $S'\into S$.\qed

\end{say}

\section{R-equivalence in families}

\begin{defn}\label{loc.field.defn}
 A {\it local field} $K$ is the quotient field
of a complete local Dedekind ring with finite residue field.
These are the finite extensions of the $p$-adic fields $\q_p$
and the Laurent series fields over finite fields $\f_q((t))$.

The valuation determined by the unique maximal ideal
 gives a locally compact metric on $K^n$.
We call this the   $p$-adic metric (though for $\f_q((t))$
the $t$-adic metric would be a more appropriate name).
If $Z$ is a variety over $K$, its $K$-points  inherit a well defined
topology from its affine pieces, called the {\it $p$-adic topology}.
$Z(K)$ is locally compact and if $Z$ is projective then
$Z(K)$ is compact.
\end{defn}

\begin{say}[Proof of (\ref{semicont.thm})]

 The assertion is local on $Y(K)$, so we may fix
a point $0\in Y(K)$ and prove upper semi continuity there.
In the process we are allowed to replace $Y$ with any other
$(0\in Y')\to (0\in Y)$ which is \'etale at $0$.

The proof proceeds in 2 steps. First we find finitely many 
easy to handle open subsets $W_i\subset X_0(K)$
such that each $W_i$ is contained in a single
R-equivalence class and the $W_i$ lift 
to nearby fibers $X_y$.
We are done if the $W_i$ are precisely the R-equivalence classes.

If not, then there are finitely many other direct R-equivalences
$$
W_i\ni u_i\sim u_j\in W_j
$$
between points of different open sets  which generate
the full R-equivalence relation. We prove that each of these
lifts to  every fiber over a suitable
 $p$-adic neighborhood of $0 \in U_{ij}\subset Y(K)$.
Thus over any point of the intersection $\cap U_{ij}$
the number of R-equivalence classes is at most as big as
over $0$.

By \cite{rcloc},
for any point $x\in X_0(K)$ there is a 
morphism $g_x:\p^1\to X_0$ over $K$  such that
\begin{enumerate}
\item $g_x(0{:}1)=x$,
\item $g_x(1{:}0)=:x'\neq x$, and
\item $g_x^*T_{X_0}(-2)$ is ample.
\end{enumerate}

Consider the space of morphisms with universal map
$$
u: \Hom(\p^1, X_0, (1{:}0)\mapsto x')\times \p^1\to X_0.
$$
By \cite[II.3.5]{rcbook}, there is a neighborhood
$$
[g_x]\in V_{x,0}\subset \Hom(\p^1, X_0, (1{:}0)\mapsto x')
$$
such that the evaluation map
$$
u^{(0{:}1)}: V_{x,0}\to X_0,\qtq{given by} v\mapsto u(v,(0{:}1))
$$
is smooth.  Smooth maps are open in the $p$-adic topology, hence
$$
W_{x,0}:= u^{(0{:}1)}(V_{x,0}(K))\subset X_0(K)
$$
is $p$-adic open for every $x$. Since $X_0(K)$ is compact,
finitely many of these cover $X_0(K)$.
Let these correspond to the points $x_1,\dots,x_d$.
Note that every point in $W_{x,0}$ is directly R-equivalent to
$x'$, thus in fact any two points in $W_{x,0}$ are R-equivalent.

After an \'etale base change we can assume that there are sections
$s_i:Y\to X$ such that $s_i(0)=x'_i$.

Look at the relative space of morphisms with universal map
$$
u: \Hom(\p^1, X, (1{:}0)\mapsto s_i(Y))\times \p^1\to X.
$$
As before, there are neighborhoods
$$
[g_x]\in V_x\subset \Hom(\p^1, X, (1{:}0)\mapsto s_i(Y))
$$
such that the evaluation map
$$
u: V_x\to X_0,\qtq{given by} v\mapsto u(v,(0{:}1))
$$
is smooth in every fiber and 
$V_x\cap \Hom(\p^1, X, (1{:}0)\mapsto x'_i)=V_{x,0}$.

Set $W_{x_i}:= u(V_{x_i}(K))\subset X(K)$.
As before, for every $y$, the intersection  $W_{x_i}\cap X_y(K)$
is a $p$-adic  open set  which is contained in a  single
R-equivalence class. 
By construction  $\cup_i W_{x_i}\supset X_0(K)$,
hence there is a $p$-adic neighborhood $0\in U_1\subset Y(K)$
such that
$$
\cup_i W_{x_i}\supset X_y(K)\qtq{for every $y\in U_1$.}
$$
This accomplishes the first part of the proof.

Let now $h:\p^1\to X_0$ be any morphism
and set $u:=h(0{:}1)$ and
$v:=h(1{:}0)$. 
After an \'etale base change we can assume that there are sections
$\sigma_u,\sigma_v:Y\to X$ such that $\sigma_u(0)=u$ and $\sigma_v(0)=v$.

Let $Z_{u,v,n}$ denote the blow up of $X\times \p^n$
along 
$$
\sigma_u(Y)\times (0{:}\cdots{:}0{:}1)\cup
\sigma_v(Y)\times (1{:}0{:}\cdots{:}0),
$$
and $E_u,E_v\subset Z_{u,v,n}$ the exceptional divisors.
As  in (\ref{main.2.pf}) there is an $n\geq 1$ 
and a genus 0 comb $C^*\subset Z_{u,v,n}$
such that $C^*$ intersects the central fibers
of $E_u,E_v$ in a single point transversally
and that $[C^*]$ is a smooth point of the relative Hilbert scheme
$\hilb(Z_{u,v,n}/Y)$. As before, this implies that
the points $\sigma_u(y),\sigma_v(y)\in X_y(K)$
are R-equivalent for every $y$ in a suitable 
$p$-adic neighborhood of $0$.

This proves the second part and thereby also (\ref{semicont.thm}).\qed
\end{say}

\section{R-equivalence over large fields}

The difference between R-equivalence and direct R-equivalence
is a frequent source of technical problems.  The aim of this section
is to prove that these two notions coincide over  large fields.

\begin{defn} \label{large.def}
A field $K$ is called {\it large} if every $K$-variety
with a smooth $K$-point has a dense set of $K$-points.

The best known examples are local fields and
infinite algebraic extensions of finite fields.
\end{defn}

\begin{thm}\label{large.fld.thm}
 Let $K$ be a large field and $X$  a smooth, projective, SRC variety over $K$
with $\dim X\geq 3$.
 For any set of points  $x_1,\dots, x_n\in X(K)$ 
the following two assertions are equivalent:
\begin{enumerate}
\item There is a smooth rational curve $C\subset X$ containing all
 the points $x_i$ such that  $T_X\otimes \o_C(-n)$ is ample.
\item Any two of the points are R-equivalent.
\end{enumerate}
\end{thm}

Proof. It is clear that (\ref{large.fld.thm}.1) implies 
(\ref{large.fld.thm}.2), the interesting part is the converse.

Let us first study the situation for two directly R-equivalent
points $x,y\in X(K)$. The direct R-equivalence is given by
some morphism $f:\p^1\to X$ with $f(0{:}1)=x, f(1{:}0)=y$.
Consider the comb  $C^*\subset X\times \p^n$ constructed in 
(\ref{SRC.noobs.thm}).
As in (\ref{main.2.pf}) we blow up the points corresponding to
$x,y$ and look at the irreducible component $H(C^*)$  of the Hilbert scheme 
of $B_{xy}(X\times \p^n)$
containing the birational transform of $C^*$.  
By general results on the Hilbert scheme
(see, for instance, \cite[37]{arko})
$[C^*]$ is a smooth $K$-point
of $H(C^*)$
and a general point of of $H(C^*)$ corresponds to a smooth rational curve
with nef normal bundle.
Since $K$ is a large field, 
   $K$-points are dense in $H(C^*)$,
hence there is a smooth rational curve 
$ A\subset B_{xy}(X\times \p^n)$ defined over $K$ with nef normal bundle.
$A$ intersects each exceptional divisor in a single point,
so $A\cong \p^1$.
Let $g:\p^1\cong A \to X$ be the projection.
Undoing the 2 blow ups twists the normal bundle by $\o(2)$,
 thus, after suitable reparametrization, 
  $g(0{:}1)=x, g(1{:}0)=y$ and $g^*T_X$ is ample.

This seems a relatively small advance, but the deformation theory 
of $g$ is much better than for $f$. 

Assume now that we have 3 points $x,y,z$ and $x,y$ and $y,z$ are
directly R-equivalent. By the above arguments, these R-equivalences
can be realized by morphisms 
$g,h:\p^1\to X$ such that $g(0{:}1)=x, g(1{:}0)=y$,
 $h(0{:}1)=y, h(1{:}0)=z$ and $g^*T_X, h^*T_X$ are ample.
The  gluing technique (cf. \cite[16]{arko}) now implies
that  there is a $K$-morphism
$f_{xz}:\p^1\to X$ such that $f_{xz}(0{:}1)=x, f_{xz}(1{:}0)=y$ and 
$f_{xz}^*T_X$ is ample. By iterating this procedure we obtain that if
two points $x,y\in X(K)$ are R-equivalent then
 there is a $K$-morphism
$f_{xy}:\p^1\to X$ such that $f_{xy}(0{:}1)=x, f_{xy}(1{:}0)=y$ 
and $f_{xy}^*T_X$ is ample.

Given our points $x_1,\dots,x_n\in X(K)$ choose another point
$y$ which is R-equivalent to them. (We have at least one rational curve in $X$
through any point, 
so there are infinitely many points in every R-equivalence class.)
There are   $K$-morphisms
$f_i:\p^1\to X$ such that $f_i(0{:}1)=x_i, f_i(1{:}0)=y$ and
 $f_i^*T_X$ is ample.

We can now proceed as in \cite[42]{arko}. Out of the morphisms $f_i$
we assemble a comb $C=C_0\cup\cdots\cup C_n$ with 
$C_i\cong \p^1$ and a morphism
$F:C\to X$ such that $F$ maps $C_0$ to $y$ and $F|_{C_i}=f_i$.
In \cite[42]{arko} we only assumed that the bundles $f_i^*T_X$ are
nef, and here they are even ample. This implies that 
deformations of $(C,F)$ give morphisms $f:\p^1\to X$
passing through all the points $x_i$ and such that $f^*T_X(-n)$ is ample.

If $\dim X\geq 3$ then we can even find such an $f$ which is an embedding
by \cite[II.3.14]{rcbook}.\qed

\section{R-equivalence on quartics}

It is conjectured that if $K$ is finitely generated
over $\q$ and $X$ is a $K$-variety then
$\ch_0(X)$ is finitely generated. This would imply that
if $X$ is also SRC then $\ch_0(X)$ is $\z + (\mbox{finite group})$.

This  suggests that under similar conditions,
$X(K)/R$ might be finite. 
The aim of the next example is to prove that this
is not so, at least over fields of transcendence
degree at least 1.  It is still possible
that if $K$ is a number field and $X$ is rationally connected
then $X(K)/R$ is finite. This is open even for surfaces.

\begin{exmp}\label{quartic.const}
 Let $Y/\q$ be any smooth variety.
By an observation  of \cite{mumf}, there is
an embedding $Y\subset \p^n$
such that $Y$ is defined by quadratic equations
$$
Y:=(q_1=\cdots=q_m=0)
$$
 and $Y$ is contained in a  smooth quartic hypersurface $(g=0)$.
Consider  the hypersurface
$$
X_{\q[t]}:= (q_1^2+\cdots+q_m^2+tg=0)\subset \p^n_{\spec \q[t]}.
$$
$X_{\q[t]}$ is a quartic hypersurface. The generic fiber 
$X_{\q(t)}$ is smooth
and the special fiber is a singular quartic
$$
X_{\q}=(q_1^2+\cdots+q_m^2=0).
$$
By our construction, the real points of $X_{\q}$
coincide with the real points of $Y$, thus
the same holds for $\q$-points:
$$
X_{\q}(\q)=Y(\q).
$$

\begin{claim} There is a surjective specialization map
$$
X_{\q(t)}(\q(t))/R\onto Y(\q)/R.
$$
\end{claim}

Proof. The specialization of a $\q(t)$-point of $X_{\q(t)}$
is a $\q$-point of $X_{\q}$, but these are all in $Y(\q)$.
Conversely, if $(a_0:\cdots:a_n)$ is a $\q$-point
of $Y(\q)$ then this also gives a point
of $X_{\q[t]}$ since both $g$ and the $q_i$
vanish on $Y$.
Thus we have a  surjective specialization map
$$
X_{\q(t)}(\q(t))/R\onto X_{\q}(\q)/R.
$$
The injection $Y\into X_{\q}$ gives a map
$$
Y(\q)/R\to X_{\q}(\q)/R 
$$
which is surjective. In order to establish our claim,
we still need to establish that two
$\q$-points of $Y$ are 
 R-equivalent in $X_{\q}$ iff they are R-equivalent in $Y$.

To see this, let $h:\p^1\to X_{\q}$ be any morphism.
$\q$-points of $\p^1$ are Zariski dense, thus
$h(\p^1)$ is contained in the Zariski closure of
$X_{\q}(\q)$. This is, however, contained in  $Y$.
Thus any chain of rational curves in $X_{\q}$ is entirely inside
$Y$.\qed
\end{exmp}

\begin{say}[Proof of (\ref{nonfg.exmp})]

With the above notation, 
let us  consider the case when $Y=E\cup Z$ is the disjoint union of
an elliptic curve over
$\q$ and another variety $Z$.
There are no non-constant maps $\p^1\to E$, thus the above construction
gives a quartic hypersurface $H:=X_{\q(t)}$ with
a surjection
$$
X_{\q(t)}(\q(t))/R\onto  E(\q)+Z(\q)/R.
$$
If we choose $E$ to have   infinitely many $\q$-points
then $X_{\q(t)}$
has infinitely many R-equivalence classes.

In order to get a unirational example, choose any 
elliptic curve over
$\q$ with infinitely many $\q$-points.
We embed $E\into \p^m$ by any complete linear system of degree
$m+1\geq 4$. Then $E$ is an intersection of quadrics.
(For $m=3$ it is a complete intersection of 2 quadrics, this case
is enough for us.) Let $x_0,\dots,x_m$ be coordinates on
$\p^m$ and $q_r=0$ the equations of $E$.

Consider now $\p^{m+s+1}$ with coordinates
$x_0,\dots,x_m, y_0,\dots,y_s$ and the equations
$q_r(x_0,\dots,x_m)$ and $x_iy_j=0$ for every $i,j$.
These equations define a scheme $Y$ with 2 irreducible components.
One is $E$ embedded in the linear space $(y_0=\cdots=y_s=0)$
and the other is the linear space $L=(x_0=\cdots=x_m=0)\cong \p^s$.
$Y=E\cup L$ is an intersection of quadrics and
it is contained in a  smooth quartic if $s\leq m$.

We can construct $X_{\q(t)}\subset \p^{m+s+1}$ as above.
$X_{\q(t)}$ has infinitely many R-equivalence classes
and it contains a linear space of dimension $s$.

Any smooth quartic of dimension at least 3 is Fano
and rationally connected. 
Unirationality of quartics containing a 2-plane
is outlined in \cite[V.5.18]{rcbook}.
This gives a degree 6 map $\p^{m+s}\map  X_{\q(t)}$.

Thus for every $n\geq 5$ we obtain a unirational quartic
of dimension $n$ with infinitely many R-equivalence classes.

The Chow group of this example does not seem to be
very interesting. A degree zero 0-cycle
 is 6-torsion and  we only get a
surjection from $\ch_0^0(H)$ to a torsion quotient
of $E(\q)$ and the latter is finite.
The relevant quotient is probably $E(\q)/2E(\q)$.
\qed
\end{say}

\medskip

\begin{say}[Proof of (\ref{deg1.zeroc.exmp})]

Let $k/\q$ be a field extension of degree $d$. 
Construct a quartic $X_d$ over $\q(t)$ using $Y=\spec k$ as in
 (\ref{quartic.const}).
$X_d$ has  a point over $k(t)$ and $\deg k(t)/\q(t)=d$.

Assume that $X_d$ has a point $P$
over a field extension  $L'/\q(t)$ of odd degree $d'$.
Specialize this point to a zero cycle 
$\sum m_i[p_i]$ of degree $d'$ on
$Y_d$.  One of the  $p_i$, say $p_1$,  has odd degree.
$p_1\not\in X_{\q}\setminus Y_{\q}$ since otherwise
we would obtain a point of the quadric
$u_1^2+\cdots+u_m^2=0$ over an odd degree extension.
Thus $p_1$  is contained
in $Y=\spec k$,
 hence its degree is divisible by $d$. Thus $d'\geq d$.
\qed
\end{say}

\begin{ack}   I   thank J.-L.\ Colliot-Th\'el\`ene 
and E.\ Szab\'o
for helpful
comments. This paper was completed during my stay
at the Isaac Newton Institute for Mathematical Sciences, Cambridge.
Partial financial support was provided by  the NSF under grant number 
DMS-9970855. 
\end{ack}

\vskip1cm

\noindent Princeton University, Princeton NJ 08544-1000

\begin{verbatim}kollar@math.princeton.edu\end{verbatim}

\end{document}